\documentclass[10pt, a4paper]{amsart}

\usepackage[utf8]{inputenc}
\usepackage[english]{babel}

\usepackage{amsthm}
\usepackage{amstext}
\usepackage{amsmath}
\usepackage{amsfonts}
\usepackage{amssymb}
\usepackage{mathtools}

\usepackage{tikz}
\usepackage{graphics}
\usepackage{wrapfig}
\usetikzlibrary{matrix,arrows,decorations.pathmorphing}
\usepackage{comment}

\usepackage[a4paper]{geometry}
\usepackage{tikz-cd}
\usepackage{epigraph}

\usepackage{hyperref}
\usepackage{cite}

\usepackage{enumerate}
\usepackage{enumitem}

\author{Matteo Tamiozzo}
\title{Special geodesics and atypical intersections}

\swapnumbers
\newtheorem{thm}[subsubsection]{Theorem}     
\theoremstyle{plain}                    
\theoremstyle{definition}               
\newtheorem{defin}[subsubsection]{Definition}
\theoremstyle{remark}                   
\newtheorem{rem}[subsubsection]{Remark}      

\newcommand{\R}{\mathbf{R}}     
\newcommand{\C}{\mathbf{C}}     
\newcommand{\Q}{\mathbf{Q}}      
\newcommand{\Z}{\mathbf{Z}}    

\newcommand{\Hp}{\mathbf{H}}

\numberwithin{equation}{subsection}

\begin{document}

\begin{abstract}
Let $C$ be a complex irreducible plane curve that is not the vanishing locus of a modular polynomial. We show that $C$ contains finitely many real algebraic curves whose projection on each coordinate axis is a union of special geodesics.
\end{abstract}

\maketitle

\tableofcontents

\section{Introduction}

\subsection{Special properties of real quadratic geodesics} In Darmon and Vonk's approach to explicit class field theory for real quadratic fields \cite{dv21}, geodesics in the upper half-plane $\Hp$ whose endpoints in $\mathbf{P}^1(\mathbf{R})$ are real quadratic conjugates play a role analogous to the one of Heegner points in the theory of complex multiplication. In view of this, one is led to wonder which properties making Heegner points ``special'' are shared by the above geodesics. Let us mention two such properties, of different flavor.
\begin{enumerate}
\item Formally analogous equidistribution statements hold true for Heegner points and for geodesics with conjugate real quadratic endpoints: see \cite{du88} - and \cite{per25} for a recent development in this direction - as well as \cite{mic04}, \cite{nord23}.
\item Heegner points in $\Hp \subset \C$ are the only algebraic points whose image via the $j$-invariant is algebraic. A suitable ``geometric'' counterpart of this statement holds true for geodesics with conjugate real quadratic endpoints, as well as for those with rational endpoints \cite[\S 3]{tam23}.
\end{enumerate}

\subsection{Detecting strongly special curves} Let $Y(1)$ be the modular curve with full level, isomorphic to $\mathbf{A}^1_\C$. For every integer $N \geq 1$, the modular polynomial $\Phi_N(T_1, T_2)$ defines an irreducible complex algebraic curve inside $Y(1)\times Y(1)\simeq \mathrm{Spec}(\C[T_1, T_2])$, image of the Hecke correspondence $Y_0(N)\rightarrow Y(1) \times Y(1)$. A curve of this form is called \emph{strongly special}. Another feature of Heegner points is that they can be used to detect strongly special curves in $Y(1)\times Y(1)$. Indeed, André proved the following result - a special case of the André--Oort conjecture.
\begin{thm}(André \cite{and98})\label{thm-and}
An irreducible complex algebraic curve in $Y(1)\times Y(1)$ that is not horizontal nor vertical is strongly special if and only if it contains infinitely many points whose coordinates are Heegner points.
\end{thm}

\subsubsection{} Let us call \emph{real quadratic geodesic} in $\Hp$ a geodesic with conjugate real quadratic endpoints; the image via $j: \Hp \rightarrow \C=Y(1)(\C)$ of a real quadratic geodesic in $\Hp$ will be called a real quadratic geodesic in $Y(1)$. In line with the above discussion, one may ask for a characterization of strongly special curves in $Y(1) \times Y(1)$ involving real quadratic geodesics rather than Heegner points. Precisely, Darmon raised\footnote{After Scanlon's talk ``What is… O-Minimality?'' at MSRI, on May 15th, 2023.} the following question.

\subsubsection{Question}\label{ssub-mainq} Let $C \subset Y(1) \times Y(1)$ be an irreducible complex algebraic curve. Assume that $C$ contains infinitely many real algebraic curves whose projection to each copy of $Y(1)$ is a union of real quadratic geodesics. Is $C$ a strongly special curve?

\subsubsection{} The aim of this note is to give a positive answer to the above question; this will be achieved in \ref{ssec-answ}. In fact, our path to the answer will lead us to a more general result, Theorem \ref{thm-main} below.

The key idea behind our approach to Question \ref{ssub-mainq} is to interpret it from an intersection-theoretic viewpoint. In the question, we are implicitly identifying the complex points of $Y(1)$ with the real points of $\mathbf{A}^2_\R$. The curve $C$ gives rise via Weil restriction to a real surface $\tilde{C}$ in $\mathbf{A}^2_\mathbf{R}\times \mathbf{A}^2_\mathbf{R}$. We are interested in real algebraic curves - whose properties will be made more precise in \ref{ssec-answ} - in $\mathbf{A}^2_\R \times \mathbf{A}^2_\R$, whose real points are contained in sets of the form $\mathcal{S} \cap \tilde{C}(\mathbf{R})$, where $\mathcal{S} \subset \mathbf{R}^2 \times \mathbf{R}^2$ is a union of products of real quadratic geodesics in $Y(1)$. Note that the intersection $\mathcal{S} \cap \tilde{C}(\R)$ is ``atypically large'': as $\tilde{C}(\mathbf{R})$ and $\mathcal{S}$ are two-dimensional subsets of $\mathbf{R}^4$, one would expect their intersection to be zero-dimensional.

This viewpoint suggests to relate the above question to the Zilber--Pink conjecture, predicting finiteness of the set of ``atypical components'' of a complex subvariety of a Shimura variety. We obtain such a relation via the base change strategy already employed for different purposes in \cite{tam23}; we then deduce a positive answer to Question \ref{ssub-mainq} from a known case of a weak version of the Zilber--Pink conjecture for $Y(1)^4$ \cite{pits16}.

\begin{rem}
Note that, if $C$ is strongly special, then $\tilde{C}$ does contain infinitely many real algebraic curves as in Question \ref{ssub-mainq}, which can be constructed intersecting $\tilde{C}\subset \mathbf{A}^2_\mathbf{R}\times \mathbf{A}^2_\mathbf{R}$ with suitable real curves $Z_M$ in the first copy of $\mathbf{A}^2_\mathbf{R}$, cf. \ref{subsubsec-infman}.
\end{rem}

\subsection{Acknowledgments} Question \ref{ssub-mainq}, which led me to the results contained in this document, was brought to my attention by Darmon after my talk at his birthday conference ``Arithmetic cycles, Modular forms, and L-functions''. I am very grateful to the organizers for giving me the opportunity to speak at the conference. I wish to sincerely thank Darmon for sharing Question \ref{ssub-mainq}, and for suggesting that it could be related to what I had discussed in my talk. It is a pleasure to dedicate this note to Darmon on the occasion of his 60th birthday, with gratitude and admiration.

\section{Real atypical intersections in $Y(1)\times Y(1)$}

\subsection{The Zilber--Pink conjecture for $Y(1)^n$} Let us quickly recall the terminology and the result from \cite[\S 7]{pits16} that we will need. Fix an integer $n \geq 1$, and let $Y(1)^n=\mathrm{Spec}(\C[T_1, \ldots, T_n])$. By a subvariety of $Y(1)^n$ we always mean an irreducible closed subset (with its reduced subscheme structure). A special subvariety of $Y(1)^n$ is an irreducible component of a closed subset defined by finitely many equations either of the form $\Phi_{N}(T_i, T_j)=0$, for $1 \leq i < j \leq n$ and $N\geq 1$, or of the form $T_i=c$ for $1 \leq i \leq n$, where $c \in \C$ is a Heegner point.

Let $V\subset Y(1)^n$ be a subvariety. An atypical subvariety of $V$ is an irreducible component $A$ of an intersection $V \cap S$, where $S$ is a special subvariety of $Y(1)^n$, such that
\begin{equation*}
\dim(A) > \dim(V)+\dim(S)-n.
\end{equation*}
We say that $A$ is strongly atypical if it is atypical and no coordinate is constant on $A \subset Y(1)^n$.
\begin{rem}\leavevmode\label{rem-strat}
\begin{enumerate}
\item Note that the above inequality is equivalent to the inequality $\mathrm{codim}(A)<\mathrm{codim}(V)+\mathrm{codim}(S)$, where codimensions are taken in $Y(1)^n$; intuitively, this tells us that $V \cap S$ is ``atypically large''.
\item Observe that $V$ is an atypical subvariety of itself if and only if it is contained in a special subvariety $S$ different from $Y(1)^n$.
\end{enumerate}
\end{rem}
\begin{thm}(Pila--Tsimerman, \cite[Theorem 7.1]{pits16})\label{thm-pits} Every subvariety $V \subset Y(1)^n$ contains finitely many maximal strongly atypical subvarieties.
\end{thm}
\begin{rem}\label{rem-zp}\leavevmode
\begin{enumerate}
\item The Zilber--Pink conjecture \cite[Conjecture 7.1]{pits16} predicts that the theorem still holds true if one replaces strongly atypical subvarieties by atypical ones; see \cite{pil22} and \cite{daor25} for background and further results on this conjecture.
\item Pila and Tsimerman's proof of the above theorem relies on $o$-minimality. A different proof, making use of differential algebra, is given in \cite[\S 5]{asl22} (cf. Theorem 5.2 and Remark 5.5).
\end{enumerate}
\end{rem}

\subsection{Special geodesics in $Y(1)$} As in \cite[Definition 3.3.5, Remark 3.3.6]{tam23}, we call \emph{special geodesic} in $\Hp$ a curve $\mathcal{S}_\mathsf{A}=\{z \in \Hp \mid \mathsf{A}z=\bar{z}\}$ for some matrix $\mathsf{A} \in \mathrm{M}_2(\Q)$ with trace zero and strictly negative determinant. If $\mathcal{S}_\mathsf{A} \subset \Hp$ is a special geodesic, we call $j(\mathcal{S}_\mathsf{A})$ a special geodesic in $Y(1)$.
\begin{rem}\label{rem-heckespecgeod}
The action on $\Hp$ of a matrix $\mathsf{B}\in \mathrm{GL}_2(\Q)$ with positive determinant sends special geodesics to special geodesics: indeed, we have $\mathsf{B}\cdot \mathcal{S}_\mathsf{A}=\mathcal{S}_{\mathsf{B}\mathsf{A}\mathsf{B}^{-1}}$.
\end{rem}

\subsubsection{}\label{subsubsec-specgeod} Special geodesics in $Y(1)$ are contained in real algebraic curves $Z_N$ defined as follows: for an integer $N \geq 1$, let $Z_N\subset \mathbf{A}^2_{\R}$ be the curve with equation $\Phi_N(X+iY, X-iY)=0$. These curves enjoy the following properties:
\begin{enumerate}[label=(\roman*)]
\item for every $N \geq 1$, the curve $Z_N$ is geometrically integral;
\item every special geodesic in $Y(1)$ is contained in $Z_N$ for some $N \geq 1$;
\item for every $N \geq 1$, the set $Z_N(\R)$ is a union of special geodesics.
\end{enumerate}
The first two properties are proved in \cite[\S 3.5.1]{tam23}. To prove the third, note that if $(x, y) \in \R^2$ belongs to $Z_N(\R)$, then the elliptic curves with $j$-invariants $x+iy$ and $x-iy$ are related by an isogeny with cyclic kernel of cardinality $N$. Writing $x+iy=j(z)$, we have $x-iy=j(-\bar{z})$, and one sees as in \cite[\S 3.5.3]{tam23} that there is a matrix $\mathsf{A} \in \mathrm{M}_2(\mathbf{Z})$ with determinant $-N$ such that $\mathsf{A}z=\bar{z}$.

\subsubsection{}\label{subsubs-deff} We will later make use of the map $f: \mathrm{Spec}(\C[X, Y]) \rightarrow \mathrm{Spec}(\C[X, Y])$ induced by the morphism of $\C$-algebras sending $(X, Y)$ to $(X+iY, X-iY)$. For $N \geq 1$, the base change to $\C$ of $Z_N$ is the preimage via $f$ of the strongly special curve with equation $\Phi_N(X, Y)=0$.

\subsection{Curves with special geodesic projections}

Let $p_1: \mathbf{A}^2_\R\times \mathbf{A}^2_\R \rightarrow \mathbf{A}^2_\R$ and $p_2: \mathbf{A}^2_\R\times \mathbf{A}^2_\R \rightarrow \mathbf{A}^2_\R$ denote the two projections. By a real algebraic curve in $\mathbf{A}^2_\R\times \mathbf{A}^2_\R$ we mean a reduced, closed subscheme of dimension one.

\begin{defin}\label{def-cspecgpr}
An infinite subset $\mathcal{Z} \subset \R^2 \times \R^2$ is a curve with special geodesic projections if it has the following properties:
\begin{enumerate}
\item there is an irreducible real algebraic curve $Z \subset \mathbf{A}^2_\R \times \mathbf{A}^2_\R$ such that $Z(\R)=\mathcal{Z}$;
\item the images $p_1(\mathcal{Z})$ and $p_2(\mathcal{Z})$ are not singletons, and are contained in a union of special geodesics.
\end{enumerate}
\end{defin}

\subsubsection{}\label{subsubs-rqwit} Let $\mathcal{Z} \subset \R^2 \times \R^2$ be an infinite set such that $p_1(\mathcal{Z})$ and $p_2(\mathcal{Z})$ are not singletons. The set $\mathcal{Z}$ is a curve with special geodesic projections if and only if there exist integers $N_1, N_2 \geq 1$ and an irreducible real algebraic curve $Z \subset Z_{N_1}\times Z_{N_2}$ such that $\mathcal{Z}=Z(\R)$. Indeed, if $\mathcal{Z}$ is of this form, then $p_i(\mathcal{Z}) \subset Z_{N_i}(\R)$ for $i\in \{1, 2\}$, and $Z_{N_i}(\R)$ is a union of special geodesics by \ref{subsubsec-specgeod}(iii). Conversely, assume that $\mathcal{Z}$ has special geodesic projections. Then, by \ref{subsubsec-specgeod}(ii), the set $\mathcal{Z}=Z(\R)$ is contained in the union of the sets $Z_{N_1}(\R)\times Z_{N_2}(\R)$, as $N_1$ and $N_2$ range over all strictly positive integers. For a fixed couple $(N_1, N_2)$, the intersection $Z \cap (Z_{N_1}\times Z_{N_2})$ is a Zariski closed subset of $Z$, hence it is either $Z$ or a finite set. The set $\mathcal{Z}$, being infinite, contains a smooth point of $Z$, hence an open neighborhood in the Euclidean topology homeomorphic to an interval. In particular, $\mathcal{Z}$ is uncountable, hence there exists a couple $(N_1, N_2)$ such that $Z \subset Z_{N_1}\times Z_{N_2}$.

The following is our main result, which can be thought of as analogous to Theorem \ref{thm-and}.

\begin{thm}\label{thm-main}
An irreducible complex algebraic curve $C \subset Y(1)\times Y(1)$ which is not horizontal nor vertical is strongly special if and only if $C(\C)$ contains infinitely many curves with special geodesic projections.
\end{thm}

\begin{rem}
Special geodesics in $Y(1)$ contain infinitely many Heegner points (cf. the proof of \cite[Proposition 3.3.7]{tam23}); however, curves with special geodesic projections may contain only finitely many points $(x_1, y_1, x_2, y_2)$ such that $x_1+iy_1$ and $x_2+iy_2$ are Heegner points. For instance, for any $\lambda \in \R^\times$ the set $\{(x_1, y_1, x_2, y_2) \in \R^4 \mid y_1=y_2=0, x_2=\lambda x_1\}$ is a curve with special geodesic projections (cf. \cite[\S 3.5.2]{tam23}); if $\lambda$ is transcendental and $x_1 \neq 0$, then $x_1$ and $\lambda x_1$ cannot be both Heegner points.

Therefore, it is not clear to us if Theorem \ref{thm-main} can be deduced from Theorem \ref{thm-and}. Instead, as announced in the introduction we will follow a different path, and deduce (one of the implications of) Theorem \ref{thm-main} from Theorem \ref{thm-pits} for $n=4$.
\end{rem}

\subsection{Proof of the main result} We will now prove Theorem \ref{thm-main}. The restriction of scalars from $\C$ to $\R$ of $Y(1)\times Y(1)=\mathrm{Spec}(\C[T_1, T_2])$ is $\mathrm{Spec}(\R[X_1, Y_1, X_2, Y_2])$, the bijection between $\R$-points of the latter and $\C$-points of the former being the map which sends $(x_1, y_1, x_2, y_2)$ to $(x_1+iy_1, x_2+iy_2)$.

\subsubsection{}\label{subsubsec-infman} Let us first show that, if $C\subset Y(1) \times Y(1)$ is a strongly special curve, then $C(\C)$ contains infinitely many curves with special geodesic projections. The curve $C$ has equation $\Phi_N(T_1, T_2)=0$ for some $N \geq 1$; let $\tilde{C}=\mathrm{Res}_{\C/\R}C \subset \mathbf{A}^2_\R \times \mathbf{A}^2_\R$. Given an integer $M \geq 1$, consider the intersection $\tilde{C} \cap (Z_M \times \mathbf{A}^2_\R)\subset \mathbf{A}^2_\R \times \mathbf{A}^2_\R$. First of all, we claim that the images via $p_1$ and $p_2$ of $\tilde{C}(\R)\cap (Z_M(\R)\times \mathbf{R}^2)$ are contained in a union of special geodesics. For $p_1$, the claim follows from \ref{subsubsec-specgeod}(iii). To prove the claim for $p_2$ note that, for a point $(x_1, y_1, x_2, y_2) \in \tilde{C}(\R)$, two elliptic curves with $j$-invariants $x_1+iy_1$ and $x_2+iy_2$ are related by an isogeny with cyclic kernel of cardinality $N$. Writing $x_1+iy_1=j(z_1)$ and $x_2+iy_2=j(z_2)$ we see that there is a matrix $\mathsf{A} \in \mathrm{M}_2(\Z)$ with positive determinant such that $z_2=\mathsf{A}z_1$. If $z_1$ belongs to a special geodesic, then so does $z_2$ by Remark \ref{rem-heckespecgeod}.

Let $Z\subset \tilde{C} \cap (Z_M \times \mathbf{A}^2_\R)$ be an irreducible component such that $\mathcal{Z}=Z(\R)$ is infinite. The set $\mathcal{Z}$ is a curve with special geodesic projections. Indeed, observe that $\tilde{C}$ is not contained in $Z_M \times \mathbf{A}^2_{\R}$ (the projection of $C \subset Y(1) \times Y(1)$ on the first factor is $Y(1)$). As the variety $Z_M \times \mathbf{A}^2_\R$ is cut out by one equation in $\mathbf{A}^2_\R \times \mathbf{A}^2_\R$, the irreducible component $Z$ has codimension one in $\tilde{C}$ \cite[Theorem 5.32]{gw20}. Finally, the restrictions of $p_1$ and $p_2$ to $\tilde{C}$ are finite, hence the maps $\mathcal{Z}\rightarrow \R^2$ induced by $p_1$ and $p_2$ have finite fibers. It follows that the sets $p_1(\mathcal{Z})$ and $p_2(\mathcal{Z})$ are not singletons. 

Varying $M$, we obtain infinitely many curves with special geodesic projections contained in $\tilde{C}(\R)=C(\C)$.

\subsubsection{Setup} Now let $C\subset Y(1) \times Y(1)$ be an irreducible complex algebraic curve whose complex points contain infinitely many curves with special geodesic projections. We will prove that $C$ is strongly special. Let $\tilde{C}=\mathrm{Res}_{\C/\R}C\subset \mathrm{Spec}(\R[X_1, Y_1, X_2, Y_2])$. We will use below the fact that, as $C$ is not horizontal nor vertical, its projection on each copy of $Y(1)$ contains a non-empty Zariski open set \cite[Theorem 10.19]{gw20}.

Let $Z\subset \mathbf{A}^2_\R \times \mathbf{A}^2_\R$ be an irreducible real algebraic curve such that $\mathcal{Z}=Z(\R)$ is a curve with special geodesic projections, and $\mathcal{Z}$ is contained in $\tilde{C}(\R)$. The intersection $Z \cap \tilde{C}$ is a closed subset of $Z$ containing infinitely many points, hence we have $Z \subset \tilde{C}$. Furthermore, by \ref{subsubs-rqwit}, there exist integers $N_1, N_2 \geq 1$ such that $Z$ is contained in $Z_{N_1, N_2}=Z_{N_1}\times Z_{N_2}$; hence, $Z$ is contained in $\tilde{C}\cap Z_{N_1, N_2}$.

\subsubsection{Base change}\label{subsubsec-bc} Let $\tilde{C}_\C\subset \mathbf{A}^2_\C\times \mathbf{A}^2_\C$ be the base change of $\tilde{C}$, and let $Z_\C$ be the base change of $Z$. Consider the map
\begin{equation*}
f^2=f \times f: \mathbf{A}^2_\C \times \mathbf{A}^2_\C \rightarrow \mathbf{A}^2_\C \times \mathbf{A}^2_\C,
\end{equation*}
where $f$ was defined in \ref{subsubs-deff}. Let us denote by $T_1, T_2, T_3$ and $T_4$ the coordinates on the target $\mathbf{A}^4_\C$ of $f^2$. Let $V=f^2(\tilde{C}_\C)$, let $A=f^2(Z_\C)$, and let $S=S_{N_1, N_2}$ be the special subvariety of $\mathbf{A}^2_\C \times \mathbf{A}^2_\C$ with equations $\Phi_{N_1}(T_1, T_2)=\Phi_{N_2}(T_3, T_4)=0$. Note that $S$ and $V$ are irreducible surfaces, and they are distinct: otherwise, we would have $\tilde{C}_\C=(Z_{N_1, N_2})\times_\R \C$, therefore $\tilde{C}(\R)=Z_{N_1}(\R)\times Z_{N_2}(\R)$. But the projection of $\tilde{C}(\R)\subset \R^2 \times \R^2$ to the first copy of $\R^2$ contains all but finitely many points.

Let us show that $A$ is a strongly atypical subvariety of $V$. First of all, the curve $Z$ is irreducible and contains infinitely many real points, hence $Z_\C$ is irreducible  \cite[\href{https://stacks.math.columbia.edu/tag/0G69}{Tag 0G69}]{stacks-project}. It follows that $A$ is irreducible, therefore it is an irreducible component of $S\cap V$ (as $S \cap V$ has dimension at most one). Finally, no coordinate is constant on $A$. Indeed, the composite of the inclusion $\R^2 \times \R^2\subset \C^2 \times \C^2$ and of $f^2$ sends $(x_1, y_1, x_2, y_2)$ to $(x_1+iy_1, x_1-iy_1, x_2+iy_2, x_2-iy_2)$. If the first or second (resp. third or fourth) coordinate was constant on $A$, then $(x_1, y_1)$ (resp. $(x_2, y_2)$) would be constant on $\mathcal{Z}$, contradicting Definition \ref{def-cspecgpr}.

\subsubsection{Strongly atypical curves in $V$} Let us now prove that $V$ contains infinitely many strongly atypical curves. As $\tilde{C}\cap Z_{N_1, N_2}$ has finitely many irreducible components, there are finitely many curves with special geodesic projections contained in a given $Z_{N_1, N_2}(\R)$ and in $\tilde{C}(\R)$. It follows that there are infinitely many mutually distinct $Z_{N_1, N_2}$ whose real points contain curves $\mathcal{Z}\subset \tilde{C}(\R)$ with special geodesic projections. Therefore, it suffices to show the following statement: take $Z\subset \tilde{C} \cap Z_{N_1, N_2}$ and $Z'\subset \tilde{C} \cap Z_{N_1', N_2'}$ such that $Z(\R)$ and $Z'(\R)$ are curves with special geodesic projections. If $Z_{N_1, N_2}$ and $Z_{N_1', N_2'}$ are distinct surfaces in $\mathbf{A}^2_\R$ then the curves $A=f^2(Z_\C)$ and $A'=f^2(Z'_\C)$ in $\mathbf{A}^4_\C$ are distinct. Suppose that $Z_{N_1}\neq Z_{N_1'}$ - if $Z_{N_2}\neq Z_{N_2'}$, the following argument applies \textit{mutatis mutandis}. Then, the complex plane curves with equation $\Phi_{N_1}(T_1, T_2)=0$ and $\Phi_{N_1'}(T_1, T_2)=0$ are different. It follows that the projection of $S_{N_1, N_2}\cap S_{N_1', N_2'}\subset \mathbf{A}^2_\C \times \mathbf{A}^2_\C$ to the first copy of $\mathbf{A}^2_\C$ has finite image. On the other hand, as no coordinate is constant on $A$, the composite of the inclusion $A \subset \mathbf{A}^4_\C$ and the projection on any coordinate has infinite image. Therefore, the curve $A \subset S_{N_1, N_2}$ cannot be contained in $S_{N_1', N_2'}$; in particular, we have $A \neq A'$.

\subsubsection{Application of Zilber--Pink}\label{ssec-appzp} By Theorem \ref{thm-pits} the strongly atypical curves in $V$ described in the previous paragraph cannot all be maximal strongly atypical subvarieties. Therefore, the surface $V$ itself must be strongly atypical. Let us deduce that $C$ is a strongly special curve.

By Remark \ref{rem-strat}, the surface $V=f^2(\tilde{C}_\C)$ is contained in a strongly special subvariety $S \subsetneq \mathbf{A}^4_\C$, which is a component of the vanishing locus of finitely many modular polynomials $\Phi_{N}(T_j, T_k)$, with $1 \leq j< k \leq 4$. To conclude, we make the following observations on the equations defining $S$.
\begin{enumerate}[label=(\roman*)]
\item There cannot be any equation $\Phi_{N}(T_j, T_k)=0$ with $j=1, k=2$. Indeed, this would imply that the points $(x_1, y_1, x_2, y_2) \in \tilde{C}(\R)$ satisfy an equation of the form $\Phi_N(x_1+iy_1, x_1-iy_1)=0$. But the projection of $C(\C)\subset \C \times \C$ on the first component contains all but finitely many points, hence it cannot be contained in $Z_N(\R)$. For the same reason, there cannot be any equation defining $S$ with $j=3, k=4$.
\item If there is an equation defining $S$ with $j=1, k=3$ then the points $(x_1, y_1, x_2, y_2) \in \tilde{C}(\R)$ satisfy an equation of the form $\Phi_N(x_1+iy_1, x_2+iy_2)=0$. Therefore, the curve $C\subset Y(1) \times Y(1)$ is contained in, hence equal to, the strongly special curve cut out by $\Phi_N$. Similarly, if there is an equation defining $S$ with $j=2, k=4$ then $C$ is a strongly special curve.
\item If there is an equation defining $S$ with $(j, k)=(2, 3)$ or $(j, k)=(1, 4)$ then the points $(x_1, y_1, x_2, y_2) \in \tilde{C}(\R)$ satisfy an equation of the form $\Phi_N(x_1-iy_1, x_2+iy_2)=0$. As $C(\C)\subset \C \times \C$ contains infinitely many points with real first coordinate, we see that the intersection between $C$ and the complex curve cut out by $\Phi_N$ is infinite. It follows that $C$ is the strongly special curve with equation $\Phi_N=0$.
\end{enumerate}

\subsubsection{Answer to Question \ref{ssub-mainq}}\label{ssec-answ} Let $\Gamma\subset \mathbf{A}^2_\R \times \mathbf{A}^2_\R$ be a real algebraic curve of pure dimension one, not necessarily irreducible. Assume that $\mathcal{S}_1=p_1(\Gamma(\R))$ and $\mathcal{S}_2=p_2(\Gamma(\R))$ are a union of special geodesics (for instance, real quadratic geodesics), and that $\Gamma(\R)$ is contained in $\tilde{C}(\R)$. Note that we could add to $\Gamma$ irreducible components with finitely many real points, all contained in $(\mathcal{S}_1\times \mathcal{S}_2) \cap \tilde{C}(\R)$, leaving $p_1(\Gamma(\R))$, $p_2(\Gamma(\R))$ unchanged and preserving the inclusion $\Gamma(\R) \subset \tilde{C}(\R)$. One way to avoid this phenomenon is to only consider curves $\Gamma$ such that $\Gamma(\R)$ has no isolated points. Intuitively, this condition tells us that $\Gamma(\R)$ really ``looks like a curve''; precisely, it implies (using \cite[Théorème 2.3.6]{bcr87}) that $\Gamma(\R)$ is a finite union of sets homeomorphic to (non-degenerate) intervals.

Assume that we have a sequence $(\Gamma_i)_{i \geq 0}$ of real algebraic curves such that $p_1(\Gamma_i(\R))$ and $p_2(\Gamma_i(\R))$ are unions of special geodesics for every $i \geq 0$. Assume that the sets $\Gamma_i(\R)$ are pairwise distinct, are contained in $\tilde{C}(\R)$, and contain no isolated points. The last property implies that the union of the real points of the irreducible components of $\Gamma_i$ with infinitely many real points is $\Gamma_i(\R)$; hence, we may discard irreducible components with finitely many real points in what follows. We will show that $C(\C)$ contains infinitely many curves with special geodesic projections. In view of Theorem \ref{thm-main}, this implies that Question \ref{ssub-mainq} has a positive answer.

As each $\Gamma_i$ has finitely many irreducible components, up to passing to a subsequence and choosing suitable irreducible components we may assume that the curves $\Gamma_i$ are irreducible and pairwise distinct. Hence, the sets $\mathcal{Z}_i=\Gamma_i(\R)$ have the following properties: they are pairwise distinct, they are contained in $\tilde{C}(\R)$, and the projections $p_1(\mathcal{Z}_i)$ and $p_2(\mathcal{Z}_i)$ are contained in a union of special geodesics. Finally, the curve $C\subset Y(1) \times Y(1)$ is not horizontal nor vertical; therefore the fibers of each projection $C \rightarrow Y(1)$ are finite. As $\mathcal{Z}_i$ is infinite, the projections $p_1(\mathcal{Z}_i)$ and $p_2(\mathcal{Z}_i)$ cannot be singletons.

\begin{rem}
Let $W$ be the restriction of scalars from $\C$ to $\R$ of $Y(1) \times Y(1)$. In the above argument, certain real subvarieties of $W$ appear naturally: for instance, in \ref{ssec-appzp}(i) we encountered the varieties $Z_N \times \mathbf{A}^2_\R\subset W$. These do contain infinitely many curves with special geodesic projections, and can be thought of as ``special real subvarieties'' of $W$; they are ruled out in our argument if we insist that we are only interested in subvarieties of $W$ coming from \emph{complex} curves in $Y(1) \times Y(1)$.
\end{rem}

\bibliographystyle{amsalpha}
\bibliography{atyp}

\providecommand{\bysame}{\leavevmode\hbox to3em{\hrulefill}\thinspace}
\providecommand{\MR}{\relax\ifhmode\unskip\space\fi MR }
\providecommand{\MRhref}[2]{%
  \href{http://www.ams.org/mathscinet-getitem?mr=#1}{#2}
}
\providecommand{\href}[2]{#2}
\begin{thebibliography}{{And}98}

\bibitem[{And}98]{and98}
Yves {Andr\'e}, \emph{{Finitude des couples d'invariants modulaires singuliers
  sur une courbe alg\'ebrique plane non modulaire}}, {J. Reine Angew. Math.}
  \textbf{505} (1998), 203--208.

\bibitem[Asl22]{asl22}
Vahagn Aslanyan, \emph{Weak modular {Zilber}-{Pink} with derivatives}, Math.
  Ann. \textbf{383} (2022), no.~1-2, 433--474.

\bibitem[BCR87]{bcr87}
J.~Bochnak, M.~Coste, and M.-F. Roy, \emph{G{\'e}om{\'e}trie alg{\'e}brique
  r{\'e}elle.}, Ergeb. Math. Grenzgeb., 3. Folge, vol.~12, Springer, 1987.

\bibitem[DO25]{daor25}
Christopher Daw and Martin Orr, \emph{Zilber-{P}ink in a product of modular
  curves assuming multiplicative degeneration}, Duke Math. J. \textbf{174}
  (2025), no.~13, 2877--2926.

\bibitem[Duk88]{du88}
William Duke, \emph{Hyperbolic distribution problems and half-integral weight
  {Maass} forms}, Invent. Math. \textbf{92} (1988), no.~1, 73--90.

\bibitem[DV21]{dv21}
Henri Darmon and Jan Vonk, \emph{Singular moduli for real quadratic fields: a
  rigid analytic approach}, Duke Math. J. \textbf{170} (2021), no.~1, 23--93.

\bibitem[GW20]{gw20}
Ulrich G{\"o}rtz and Torsten Wedhorn, \emph{Algebraic geometry {I}. {Schemes}.
  {With} examples and exercises}, 2nd ed., Springer, 2020.

\bibitem[Mic04]{mic04}
Philippe Michel, \emph{The subconvexity problem for {Rankin}-{Selberg}
  {{\(L\)}}-functions and equidistribution of {Heegner} points}, Ann. Math. (2)
  \textbf{160} (2004), no.~1, 185--236.

\bibitem[Nor23]{nord23}
Asbj{\o}rn~Christian Nordentoft, \emph{Concentration of closed geodesics in the
  homology of modular curves}, Forum Math. Sigma \textbf{11} (2023), 38.

\bibitem[Pil22]{pil22}
Jonathan Pila, \emph{Point-counting and the {Zilber}-{Pink} conjecture}, Camb.
  Tracts Math., vol. 228, Cambridge: Cambridge University Press, 2022.

\bibitem[PP25]{per25}
Patricio P{\'e}rez-Pi{\~n}a, \emph{{{\(p\)}}-adic equidistribution of modular
  geodesics and of {CM} points on {Shimura} curves}, J. Number Theory
  \textbf{271} (2025), 259--282.

\bibitem[PT16]{pits16}
Jonathan Pila and Jacob Tsimerman, \emph{Ax-{S}chanuel for the
  {{\(j\)}}-function}, Duke Math. J. \textbf{165} (2016), no.~13, 2587--2605.

\bibitem[{Sta}26]{stacks-project}
The {Stacks project authors}, \emph{The stacks project},
  \url{https://stacks.math.columbia.edu}, 2026.

\bibitem[Tam23]{tam23}
Matteo Tamiozzo, \emph{Special curves in modular surfaces}, Can. Math. Bull.
  \textbf{66} (2023), no.~1, 1--18.

\end{thebibliography}

\end{document}